\documentclass[12pt]{article}
\usepackage{amsfonts,amsmath,amscd}
\usepackage[a4paper,left=3cm,right=2cm,top=3cm,bottom=2cm]{geometry}
\usepackage{hyperref}
\usepackage{lmodern}
\usepackage[utf8x]{inputenx} %
\usepackage[LGR, T1]{fontenc}
\usepackage{textcomp} 
\usepackage{graphicx}
\usepackage{wrapfig}

\newtheorem{remark}{Remark}[section]

\newtheorem{exam}{Example}[section]
\newtheorem{teorema}{Theorem}[section]

\newcommand{\n}{\noindent}

\begin{document}

\title{A family of MCF solutions for the Heisenberg Group}

\author{
\textbf{Benedito Leandro}
\\
{\small\it Universidade Federal de Jata\'i }\\
{\small\it BR 364, km 195, 3800, 75801-615, CIEXA,
	Jata\'i, GO, Brazil. }
\\
{\small\it e-mail:  bleandroneto@gmail.com}
  \\
\textbf{Adriana Araujo Cintra}
\\
{\small\it Universidade Federal de Jata\'i }\\
{\small\it BR 364, km 195, 3800, 75801-615, CIEXA,
	Jata\'i, GO, Brazil. }
\\
{\small\it e-mail:  adriana.araujo.cintra@gmail.com}
  \\
  \textbf{Hiuri Fellipe dos Santos Reis}
\\
{\small\it Instituto Federal de Goiás}\\
{\small\it Rua Formosa, Lot. Santana, 76400-000,
	Urua\c{c}u, GO, Brazil. }
\\
{\small\it e-mail:  hiuri.reis@ifg.edu.br}
  \\
}

\maketitle
\thispagestyle{empty}

\markboth{abstract}{abstract}
\addcontentsline{toc}{chapter}{abstract}

\begin{abstract}

     \noindent
The aim of this paper is to investigate the mean curvature flow soliton solutions on the Heisenberg group $\mathcal{H}$ when the initial data is a ruled surface by straight lines. We give a family of those solutions which are generated by $\mathfrak{Iso}_{0}(\mathcal{H})$ (the isometries of $\mathcal{H}$ for which the
origin is a fix point). We conclude that the function which describe the motion of these surfaces under MCF, is always a linear affine function. As an application we proof that the Grim Reaper solution evolves from a ruled surface in $\mathcal{H}$. We also provide other examples.
         \end{abstract}

\noindent 2010 Mathematics Subject Classification: 53C44, 35R03, 14J26. \\
Keywords: Mean curvature flow; Ruled surface; Heisenberg group.

\section{Introduction and main statement}

The mean curvature flow (MCF) is a geometric evolution equation. In short, MFC is a way to let submanifolds evolve in a given manifold over time to minimize its volume. Consider a map 
$$\Phi:M^{n}\rightarrow (N^{n+1},\,\bar{g}),$$
of a given manifold $M$ on a smooth Riemannian manifold $N$ with Riemannian metric $\bar{g}$. The induced metric $g:=\{g_{ij}\}$ on $M$ is 
$$g_{ij}=\left\langle\frac{\partial\Phi}{\partial x_{i}},\,\frac{\partial\Phi}{\partial x_{j}}\right\rangle_{\bar{g}}=\bar{g}_{\alpha\beta}\frac{\partial\Phi^{\alpha}}{\partial x_{j}}\frac{\partial\Phi^{\beta}}{\partial x_{j}},$$
where $1\leq i,\,j\leq n$ and $1\leq\alpha,\,\beta\leq n+1$. Thus, the second fundamental form $h:=\{h_{ij}\}$ for this hypersurface is
\begin{eqnarray}\label{secondfund}
h_{ij}N^{\alpha}=\frac{\partial^{2}\Phi}{\partial x_{i}\partial x_{j}}-\Gamma^{k}_{ij}\frac{\partial\Phi^{\alpha}}{\partial x_{k}}+\bar{\Gamma}^{\alpha}_{\beta\gamma}\frac{\partial\Phi^{\beta}}{\partial x_{i}}\frac{\partial\Phi^{\gamma}}{\partial x_{j}},
\end{eqnarray}
where $N^{\alpha}$, $\Gamma^{k}_{ij}$ and $\bar{\Gamma}^{k}_{ij}$ are, respectively, the component of the normal vector, and the Christoffel symbols for $g$ and $\bar{g}$ (cf. \cite{Huisken,Hungerbuhler,Saez}). 

Since we know the induced metric and the second fundamental form, we can control all the geometry of this hypersurface. In particular, the mean curvature vector is given by $\Vec{H}=HN=g^{ij}\big(h_{ij}N\big)$. Therefore, from \eqref{secondfund} we have that $$\Vec{H}=\Delta\Phi,$$
where $\Delta$ is the Laplace-Beltrami operator which depends on $g$ and $\bar{g}$.  

A map $\Phi^{t}:M^{n}\times I\rightarrow N^{n+1}$, $I=[0,\,T)$,  which is a family of hypersurfaces, is a solution of MCF if $$\partial_{t}\Phi^{t}(p,\,t)=\Vec{H}^{t}(p,\,t),$$
in which $p\in M$ and $t\in I$. So, the mean curvature flow behavs pretty much like a quasilinear second order system of parabolic differential equations, similar to the ordinary heat equation. However, important differences arise on the formation of singularities (cf. \cite{Colding}).

MCF has been studied in material science to model things
such as cells, grain boundaries in annealing pure metal, bubble growth, and image processing (cf. \cite{Colding, Hungerbuhler} and the references therein). This flow and others have been extensively used to model several physical phenomena, e.g., Huisken and Ilmanen \cite{Huisken1} used the inverse MFC to proof the Riemannian Penrose inequality. However, there are very few results on MFC in non-Euclidean spaces (cf. \cite{Borisenko,Ferrari, Hiuri}).

Here, we propose to look for solutions for the mean curvature flow on the Heisenberg group $\mathcal{H}$. This group was chosen because it is very well behaved and the connections associated with the standard metric (cf. Section \ref{preliminar} for a overview of $\mathcal{H}$) for the Heisenberg group are simple, thus the derivatives are more easily calculated. Moreover, this space was already vastly revised. Werner Heisenberg introduced this group as a new approach to study quantum mechanics.

 Our intention here is to establish a simple method, able to provide explicit examples of solutions for the MCF. To do so, suppose that $\Phi^{t}:U\times I\rightarrow (\mathbb{R}^{3},\,\langle\cdot,\,\cdot\rangle)$, where $U\subset\mathbb{R}^{2}$ and $\langle\cdot,\,\cdot\rangle$ the Heisenberg metric, is a solution for the MCF with initial data $\Phi^{0}:=\Phi$, i.e.,
\begin{eqnarray}\label{3}
	\partial_{t}\Phi^{t}=\vec{H}^{t}.
\end{eqnarray}
 We say that $\Phi^{t}$ is a soliton (cf. \cite{Hungerbuhler}) if there exists a Killing vector field $X$ with flow  $\Psi^{t}$ such that
\begin{eqnarray}\label{5}
\Phi^{t}(u,\,v)=\Psi^{t}\left(\Phi(u,\,v)\right).
\end{eqnarray}
Solitons are relevant because they represent a class of solutions with very special properties, e.g., they appear as blow-ups of singularities of the MCF.

 In the Heisenberg group (as in every homogeneous manifold), usually a MCF soliton is defined as a surface that, under MCF, moves by the action of the flow of a Killing vector field. In this paper we accept this natural definition, and the solitons so defined are the analog of rotating and translating solitons in the Euclidean Space. In this paper, we find new and interesting simple families of MCF solitons in the Heisenberg group which are ruled surfaces with vertical or horizontal lines
as rulings, and, at the same time, we describe the motion of these surfaces under MCF. To find them we make use of a good knowledge of the geometry of the Heisenberg group to reduce the problem for these surfaces to an O.D.E. 

Let the initial data $\Phi^{0}$ be {\it a ruled surface} with parametrization
\begin{eqnarray}\label{1}
\Phi(u,\,v)=\alpha(u)+v\beta(u),
\end{eqnarray}
where $\alpha$ is a differentiable (but not necessarily regular)  base curve and $\beta$ is a director vector field along $\alpha$ which vanishes nowhere (cf. \cite{Kuhnel,O'neill}). Ruled surfaces are a subject wildly explored in geometry. Some important classifications of ruled surfaces were already made in Euclidean and non-Euclidean spaces (cf. \cite{Dillen, Figueroa,Shin}).

Before proceeding, we recommend the reader to see Section \ref{preliminar} to acknowledge the notation. Without further ado, we state our main results.

\begin{teorema}\label{t1}
The soliton solution for MCF on the Heisenberg group, where the $\Psi^{t}$ in definition \eqref{5} belongs to $\mathfrak{Iso}_{0}(\mathcal{H})$ for every $t$ with vertical straight lines as rulings for initial
data, is a ruled surface given by
\begin{enumerate}
	\item $\Phi^{t}(u,\,v)=\left(u\cos(At)-f(u)\sin(At),\,u\sin(At)+f(u)\cos(At),\,v\right)$,\\ where $f'(u)=\tan[A(u^{2}+f^{2}(u))+B]$.
	\
	
Moreover, we provide solutions generated by $\Psi^{t}\notin\mathfrak{Iso}_{0}(\mathcal{H})$.	
	\item $\Phi^{t}(u,\,v)=\left(At+u,\,\dfrac{\pm1}{2A}\arctan\left(\sqrt{De^{4Au}-1}\right),\,\dfrac{\pm t}{4}\arctan\left(\sqrt{De^{4Au}-1}\right)+v\right)$,

	\item $\Phi^{t}(u,\,v)=\left(u,\,\dfrac{-1}{2A}\log\left[C\sin(2Au)-B\cos(2Au)\right]+At,\,-\dfrac{At}{2}u+v\right)$,
\end{enumerate}
where $A\neq0,\,D>0,\,B,\,C\in\mathbb{R}$. 
\end{teorema}

\begin{figure}[h]
	\centering
	\includegraphics[scale=0.5]{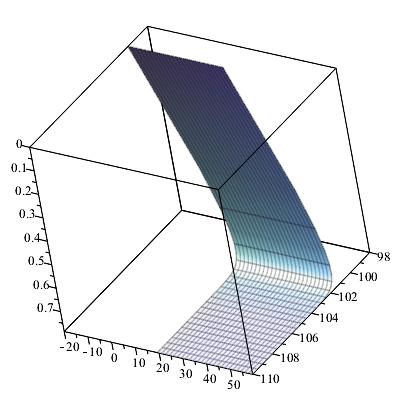}
	\includegraphics[scale=0.5]{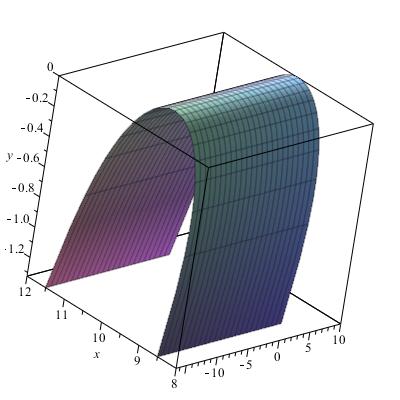}
	\caption{Soliton solution-Theorem \ref{t1}-2 /  Grim Reaper solution-Theorem \ref{t1}-3.}
\end{figure}

Horizontal straight lines are another kind of geodesic for $\mathcal{H}$ (cf. Section \ref{preliminar}). This motivates the search for this type of solitons in $\mathcal{H}$.

\begin{teorema}\label{t2}
The soliton solution for MCF on the Heisenberg group, where the $\Psi^{t}$ in definition \eqref{5} belongs to $\mathfrak{Iso}_{0}(\mathcal{H})$ for every $t$ with horizontal straight lines intersecting the $z$-axis as rulings for initial
data, is a ruled surface given by
\begin{enumerate}
    \item  $\Phi^{t}(u,\,v)=\Big(v\big(C\cos(At)-B\sin(At)\big)-f(u)\big(C\sin(At)+B\cos(At)\big),\,\\v\big(C\sin(At)+B\cos(At)\big)+f(u)\big(C\cos(At)-B\sin(At)\big),\,g(u)\Big),$
 where $f$ and $g$ are smooth functions satisfying the following equation
\begin{eqnarray*}
  &&2vAg'(u)\left\{f'(u)^{2}[4+f(u)^{2}]+[2g'(u)-vf'(u)]^{2}\right\}\nonumber\\
  &&=\left[4+f(u)^{2}\right][f''(u)g'(u)-f'(u)g''(u)].
\end{eqnarray*}
	\

Moreover, we provide solutions generated by $\Psi^{t}\notin\mathfrak{Iso}_{0}(\mathcal{H})$.

    \item  $\Phi^{t}(u,\,v)=\left(At+Cv-Bf(u),\,Cf(u)+Bv,\,\dfrac{At}{2}(Cf(u)+Bv)+g(u)\right),$
such that $f$ and $g$ are smooth functions satisfying the following equation
 \begin{eqnarray*}
	&-&A\left\{f'(u)^{2}[4+f(u)^{2}]+[2g'(u)-vf'(u)]^{2}\right\}\left[\frac{C}{2}f(u)^{2}+Bvf(u)+2Bg(u)\right]'\nonumber\\
	&=&\left[4+f(u)^{2}\right][f''(u)g'(u)-f'(u)g''(u)];
\end{eqnarray*}

\item $\Phi^{t}(u,\,v)=\left(Cv-Bf(u),\,Cf(u)+Bv,\,At+g(u)\right),$
in which $f$ and $g$ are smooth functions satisfying the following equation
 \begin{eqnarray*}
	&-&2Af'(u)\left\{f'(u)^{2}[4+f(u)^{2}]+[2g'(u)-vf'(u)]^{2}\right\}\nonumber\\
	&=&\left[4+f(u)^{2}\right][f''(u)g'(u)-f'(u)g''(u)],
\end{eqnarray*}
\end{enumerate}
where $A\neq0,\,B,\,C\in\mathbb{R}$.  
\end{teorema}

\begin{remark}
	
	\
	
	\begin{itemize}
		\item[i)] 	We can see that if either $B=C=0$, or $f$ or $g$ are constant functions, then we have trivial solutions for the MCF given by Theorem \ref{t2}. Also, if $f(u)=B g(u)+C$ we have trivial solutions.
		
		\item[ii)] In  Theorem \ref{t2}-2, if  $\frac{C}{2}f(u)^{2}+Bvf(u)+2Bg(u)=D$, where $D$ is constant, we get trivial solutions. 
		
		\item[iii)]We highlight that minimal ruled surfaces in the Heisenberg group were already classified in \cite{Shin}. They proved that parts of planes, helicoids and hyperbolic paraboloids are the only minimal surfaces ruled by geodesics in the three-dimensional Riemannian Heisenberg group. We pointed out that the quoted classification corresponds to a different definition of ruled surface.

		\item[iv)]From Theorems \ref{t1} and \ref{t2}, we can see that	any soliton solution for MCF generated by $\mathfrak{Iso}(\mathcal{H})$, in which the initial data is a ruled surface by straight lines, has a linear affine time dependent function.
	\end{itemize}
\end{remark}

In what follows, we give some examples for Theorem \ref{t2}. Even if the degree of freedom suggests that we can build solutions as much as we want, we must keep in mind that solve an ODE explicitly is not a simple task.

\begin{exam}
	In Theorem \ref{t2}-1 considering  $f(u)=u$, after make $2g'(u)-v=y(u)$ we have:
	\begin{eqnarray}\label{ddd222}
y(u)=\exp\left\{-2v\int\frac{y^{2}(u)}{4+u^{2}}du-2vu\right\}-v+E;\quad E\in\mathbb{R}.
	\end{eqnarray}
	Thus, 
	\begin{eqnarray*}
\Phi^{t}(u,\,v)=\Big(v\big(C\cos(At)-B\sin(At)\big)-u\big(C\sin(At)+B\cos(At)\big),\,\\v\big(C\sin(At)+B\cos(At)\big)+u\big(C\cos(At)-B\sin(At)\big),\,g(u)\Big),
	\end{eqnarray*}
	where $g(u)$ is implicitly given by  \eqref{ddd222}, is a soliton solution for the MCF on $\mathcal{H}$.
\end{exam}

\begin{exam}
In Theorem \ref{t2}-2 considering $C=1$, $B=0$ and $f(u)=u$, after make $2g'(u)-v=y(u)$ we have:
\begin{eqnarray}\label{ddd22}
    y'(u)=2Au\left[\frac{y(u)^{2}}{4+u^{2}}+1\right].
\end{eqnarray}
Thus, 
\begin{eqnarray*}
    \Phi^{t}(u,\,v)=\left(A t+v,\,u,\,\frac{1}{2}A tu+g(u)\right),
\end{eqnarray*}
where $g(u)$ is given by the Ricatti equation \eqref{ddd22} is a soliton solution for the MCF on the Heisenberg group.
\end{exam}

\begin{exam}
In Theorem \ref{t2}-3 consider $f(u)=u$. Then, making a change of variable $2g'(u)-v=y(u)$ we have:
\begin{eqnarray}\label{ddd2}
    y'(u)=4A\left[\frac{y(u)^{2}}{4+u^{2}}+1\right].
\end{eqnarray}
Thus, 
\begin{eqnarray*}
    \Phi^{t}(u,\,v)=\left(Cv-Bu,\,Cu+Bv,\,At+g(u)\right),
\end{eqnarray*}
where $g(u)$ is given by the Ricatti equation ODE \eqref{ddd2} is a soliton solution for the MCF on the Heisenberg group.
\end{exam}

% Demonstrações	% Demonstrações	% Demonstrações	% Demonstrações	% Demonstrações	% Demonstrações	% Demonstrações	% Demonstrações	
% Demonstrações	% Demonstrações	% Demonstrações	% Demonstrações	% Demonstrações	% Demonstrações	% Demonstrações	% Demonstrações	
% Demonstrações	% Demonstrações	% Demonstrações	% Demonstrações	% Demonstrações	% Demonstrações	% Demonstrações	% Demonstrações	
% Demonstrações	% Demonstrações	% Demonstrações	% Demonstrações	% Demonstrações	% Demonstrações	% Demonstrações	% Demonstrações	

\section{Preliminar}\label{preliminar}

The Heisenberg group $\mathcal{H}$ is $\mathbb{R}^{3}$ endowed with the left-invariant metric 
\begin{eqnarray*}
\langle\cdot,\cdot\rangle=dx^{2}+dy^{2}+\left(\frac{1}{2}ydx-\frac{1}{2}xdy+dz\right)^{2}.
\end{eqnarray*}
Its Lie algebra, in terms of the canonical basis $\{e_{1},\,e_{2},\,e_{3}\}$ of $\mathbb{R}^{3}$, is given by
\begin{eqnarray*}
[e_{1},\,e_{2}]=e_{3}\quad\mbox{and}\quad [e_{i},\,e_{3}]=0,
\end{eqnarray*}
for $i\in\{1,\,2,\,3\}$.
Using this frame, we have that an orthonormal basis of left-invariant vector fields is
\begin{eqnarray*}
E_{1}=\partial_{x}-\frac{y}{2}\partial_{z},\, E_{2}=\partial_{y}+\frac{x}{2}\partial_{z}\quad\mbox{and}\quad E_{3}=\partial_{z}.
\end{eqnarray*}
Of which
\begin{eqnarray}\label{base}
\partial_{x}=E_{1}+\frac{y}{2}E_{3},\, \partial_{y}=E_{2}-\frac{x}{2}E_{3}\quad\mbox{and}\quad\partial_{z}= E_{3}.
\end{eqnarray}

Letting $\nabla$ be the Levi-Civita connection on $\mathcal{H}$, we have
\begin{eqnarray*}
&&\nabla_{E_{i}}E_{i}=0\quad\mbox{for}\quad i=1,\,2,\,3,\quad\nabla_{E_{1}}E_{2}=\frac{1}{2}E_{3}=-\nabla_{E_{2}}E_{1},\nonumber\\
&&\quad\nabla_{E_{1}}E_{3}=-\frac{1}{2}E_{2}=\nabla_{E_{3}}E_{1}\quad\mbox{and}\quad\quad\nabla_{E_{2}}E_{3}=\frac{1}{2}E_{1}=\nabla_{E_{3}}E_{2}.
\end{eqnarray*}
Moreover,
\begin{eqnarray}\label{basevetorial}
E_{1}&=&E_{2}\times E_{3},\nonumber\\
E_{2}&=&-E_{1}\times E_{3},\nonumber\\
E_{3}&=&E_{1}\times E_{2}.
\end{eqnarray}

Now, we are ready to show the geodesics of $\mathcal{H}$. They are, basically, helices or straight lines (horizontal or vertical). We sum up this in the next result.

\begin{teorema}\cite{Marenich}
Geodesic lines issuing from the origin in $\mathcal{H}$ satisfy
the following equations
\begin{eqnarray*}
 \left\{
\begin{array}{lll}
x(u)=\frac{A}{2C}\left[\sin(2Cu+B)-\sin(B)\right], \\ \\

y(u)=\frac{A}{2C}\left[-\cos(2Cu+B)+\cos(B)\right],\\ \\

z(u)= \frac{1+C^{2}}{2C}u-\frac{1-C^{2}}{4C^{2}}\sin(2Cu),
\end{array}
\right.
\end{eqnarray*}
where $A,\,B,\,C\neq0\in\mathbb{R}.$ It is not difficult to see that when $A=0$ and $C=1$ we have straight vertical lines. And for $C=0$ the geodesics are ``horizontal'' lines and satisfy
\begin{eqnarray}\label{geohorizontal}
 \left\{
\begin{array}{lll}
x(u)=Au, \\ 

y(u)=Bu,\\ 

z(u)= 0.
\end{array}
\right.
\end{eqnarray}
\end{teorema}

 The group of isometries of the Heisenberg group  $\mathfrak{Iso}(\mathcal{H})$  was already determined in the literature (cf. \cite{Marenich, Shin}). The base of the Lie algebra of $\mathfrak{Iso}(\mathcal{H})$ is given by 
 \begin{eqnarray}\label{2}
 X_{1}=y\partial_{x}-x\partial_{y},\,X_{2}=\partial_{x}+\frac{y}{2}\partial_{z},\,X_{3}=\partial_{y}-\frac{x}{2}\partial_{z},\,X_{4}=\partial_{z}.
 \end{eqnarray}
 
 Moreover, if $\Psi^{t}$ a 1-parameter family of isometries of $\mathcal{H}$ then 
 \begin{eqnarray*}
  \left\{
 \begin{array}{lll}
 \frac{d}{dt}\Psi^{t}(p)=X(\Psi^{t}(p))\\ 
 
 \Psi^{0}(p)=p=(x,\,y,\,z),
 \end{array}
 \right.
 \end{eqnarray*}
 where $X$ is any vector field generated by \eqref{2}.
  Thus, we have the following isometries generated by the elements of \eqref{2}:
 \begin{eqnarray}\label{4}
 \left\{
 \begin{array}{lll}
 \Psi^{t}_{1}(x,\,y,\,z)=(x\cos(t)-y\sin(t),\,x\sin(t)+y\cos(t),\,z),\quad\mbox{for}\quad X_{1}; \\
 
  \Psi^{t}_{2}(x,\,y,\,z)=(t+x,\,y,\,\frac{1}{2}ty+z),\quad\mbox{for}\quad X_{2};\\
  
  \Psi^{t}_{3}(x,\,y,\,z)=(x,\,t+y,\,-\frac{1}{2}tx+z),\quad\mbox{for}\quad X_{3};\\
  
   \Psi^{t}_{4}(x,\,y,\,z)=(x,\,y,\,t+z),\quad\mbox{for}\quad X_{4}.\\
 \end{array}
 \right.
 \end{eqnarray}
 The general 1-parameter of family of isometries for a given vector field $X$ is the composition of all the above isometries. The set $\mathfrak{Iso}_{0}(\mathcal{H})$ is for those 1-parameter family of isometries for which the origin is a fix point, i.e., generated by $\Psi_{1}^{t}$.
 
 To establish our notation, the mean curvature is given by
 \begin{eqnarray}\label{meancurva}
 H=\frac{1}{2}\frac{lG-2mF+En}{EG-F^{2}},
 \end{eqnarray}
 where $E,\,F,\,G$ and $l,\,m,\,n$ stand for the coefficients of the first and second fundamental forms, respectively.
 
 Now, we are ready to proof the main results.

\section{Proof of the Main Results}

It is worth noting that the proofs of Theorem \ref{t1} and Theorem \ref{t2} seem to be long, however, the cases are recursive. That said, we divided the proof in cases, so when the reader makes the computation for the first case of each theorem, the remain cases follow easily.

\n \textbf{Proof of Theorem \ref{t1}:}

{\bf First Case:} We begin with the second item of this theorem since we believe this might be more convinient. Therefore, we consider $X_{2}$ as an element of \eqref{2}, which generates the second isometry $\Psi_{2}$ of \eqref{4}. 

Consider ruled surfaces with rulings vertical straight lines in $\mathcal{H}$ of the form (cf. \cite{Figueroa}):
\begin{eqnarray}\label{6}
\Phi(u,\,v)= (u,\,f(u),\,v),
\end{eqnarray}
where $f$ is a smooth function.
Therefore, from \eqref{5}, \eqref{4} and \eqref{6} we assume that
\begin{eqnarray}\label{7}
\Phi^{t}(u,\,v)= \left(\varepsilon(t)+u,\,f(u),\,\frac{1}{2}\varepsilon(t)f(u)+v\right)
\end{eqnarray}
is a soliton solution for MCF in $\mathcal{H}$. Here, $\varepsilon(t)$ is a smooth function such that $\varepsilon(0)=0$.

We must determinate $f(u)$ and $\varepsilon(t)$ such that they satisfy \eqref{3}. Hence, from \eqref{base} and \eqref{7} we have
\begin{eqnarray}\label{nha}
\partial_{t}\Phi^{t}&=&\left(\varepsilon'(t),\,0,\,\frac{1}{2}\varepsilon'(t)f(u)\right)\nonumber\\
&=& \varepsilon'(t)[E_{1}+f(u)E_{3}].
\end{eqnarray}

Now, we need to obtain the mean curvature $H^{t}$ and the normal vector $N^{t}$ of $\Phi^{t}$. Note that, from \eqref{base} and \eqref{7}, the first derivatives can be calculated: 
\begin{eqnarray*}
\Phi^{t}_{u}&=& E_{1}+f'(u)E_{2}+\frac{1}{2}[f(u)-uf'(u)]E_{3}\quad\mbox{and}\nonumber\\
\Phi^{t}_{v}&=&E_{3}.
\end{eqnarray*}
Thus, from \eqref{basevetorial} and the above identities, we can infer that
\begin{eqnarray}\label{nhanha}
N^{t}=\frac{\Phi^{t}_{u}\times\Phi^{t}_{v}}{\|\Phi^{t}_{u}\times\Phi^{t}_{v}\|}=\frac{f'(u)E_{1}-E_{2}}{\sqrt{1+f'(u)^{2}}}.
\end{eqnarray}

From \eqref{nha} and \eqref{nhanha}, the left-hand side of \eqref{3} has the simple form
\begin{eqnarray}\label{sup}
\langle\partial_{t}\Phi^{t},\,N^{t}\rangle=\frac{f'(u)\varepsilon'(t)}{\sqrt{1+f'(u)^{2}}}.
\end{eqnarray}

To get the mean curvature we need to prove the coefficients of the first and second fundamental forms (cf. \cite{Figueroa}). A straightforward computation ensures that  
\begin{eqnarray*}
E^{t}&=&\langle \Phi^{t}_{u},\,\Phi^{t}_{u}\rangle=1+f'(u)^{2}+\frac{1}{4}[f(u)-uf'(u)]^{2},\nonumber\\
F^{t}&=&\langle \Phi^{t}_{u},\,\Phi^{t}_{v}\rangle=\frac{1}{4}[f(u)-uf'(u)],\nonumber\\
G^{t}&=&\langle \Phi^{t}_{v},\,\Phi^{t}_{v}\rangle=1.
\end{eqnarray*}
Furthermore, 
\begin{eqnarray*}
l^{t}&=&-\langle \Phi^{t}_{u},\,\nabla_{\Phi^{t}_{u}}N^{t}\rangle=\frac{1}{2\sqrt{1+f'(u)^{2}}}[(f(u)-uf'(u)][1+f'(u)^{2}]-2f''(u),\nonumber\\\\
m^{t}&=&\langle\nabla_{\Phi^{t}_{v}}\Phi^{t}_{u},\,N^{t}\rangle=\frac{1}{2}\sqrt{1+f'(u)^{2}}\quad\mbox{and}\quad
n^{t}=\langle\nabla_{\Phi^{t}_{v}}\Phi^{t}_{v},\,N^{t}\rangle=0.
\end{eqnarray*}
So, from \eqref{meancurva} we get
\begin{eqnarray}\label{dan}
H^{t}=\frac{-f''(u)}{2\sqrt{(1+f'(u)^{2})^{3}}}.
\end{eqnarray}

Finally, from \eqref{3}, \eqref{sup} and \eqref{dan} we arrive with the ODE equation
\begin{eqnarray}\label{eq1}
f''(u)+2\varepsilon'(t)f'(u)[1+f'(u)^{2}]=0.
\end{eqnarray}
Hence, $\varepsilon(t)=At$, $A\neq0$. Moreover, \eqref{eq1} can be reduced to 
\begin{eqnarray*}
y'(u)+2Ay(u)[1+y(u)^{2}]=0,
\end{eqnarray*}
where $f'(u)=y(u)$. The solution for this ODE, is
\begin{eqnarray*}
f(u)=\frac{\pm1}{2A}\arctan\left(\sqrt{De^{4Au}-1}\right)+B,
\end{eqnarray*}
where $D:=D(A)$ is a positive constant which depends on $A$, and $B\in\mathbb{R}$.

With the first item of this theorem already proved we can reduce the cases to the following basic fact: the coefficients for the first and second fundamental forms are the same for the first, second and third items of Theorem \ref{t1}. 

{\bf Second Case:} For the third item of this theorem, considering the isometry generated by $X_{3}$, we obtain 
\begin{eqnarray}
\Phi^{t}(u,\,v)=\left(u,\,f(u)+\varepsilon(t),\,-\frac{1}{2}u\varepsilon(t)+v\right),
\end{eqnarray}
where $\varepsilon(t)$ is a smooth function such that $\varepsilon(0)=0$. Thus, an easy computation proves that
\begin{eqnarray}\label{dandan}
\partial_{t}\Phi^{t}=\varepsilon'(t)\left(E_{2}-uE_{3}\right).
\end{eqnarray}
Since the normal is also given by \eqref{nhanha}, from \eqref{dandan} we get
\begin{eqnarray*}
	\langle\partial_{t}\Phi^{t},\,N^{t}\rangle=\frac{-\varepsilon'(t)}{\sqrt{1+f'(u)^{2}}}.
\end{eqnarray*}
Combining the above equation with \eqref{dan}, the result is the following ODE:
\begin{eqnarray}\label{danada}
f''(u)-2\varepsilon'(t)[1+f'(u)^{2}]=0.
\end{eqnarray}
Then, we have $\varepsilon(t)=At$, where $A\neq0$, $C,\,B\in\mathbb{R}$, and 
\begin{eqnarray}
f(u)=\frac{-1}{2A}\log\left[C\sin(2At)-B\cos(2At)\right],
\end{eqnarray}
 such that $C^{2}+B^{2}\neq0$.

{\bf Third Case:} We now prove the first item of this theorem. From  \eqref{4} and \eqref{6} we have
\begin{eqnarray}\label{spf1}
\Phi^{t}(u,\,v)=(u\cos\varepsilon(t)-f(u)\sin\varepsilon(t),\,u\sin\varepsilon(t)+f(u)\cos\varepsilon(t),\,v).
\end{eqnarray}
The normal vector $N^{t}$ is given by
\begin{eqnarray}
N^{t}=\frac{1}{\sqrt{1+f'(u)^{2}}}[(\sin\varepsilon(t)+f'(u)\cos\varepsilon(t))E_{1}-(\cos\varepsilon(t)-f'(u)\sin\varepsilon(t))E_{2}].
\end{eqnarray}
Moreover, from \eqref{spf1} we obtain
\begin{eqnarray*}
\partial_{t}\Phi^{t}(u,\,v)&=&[\cos\varepsilon(t)-f'(u)\sin\varepsilon(t)]E_{1}+[\sin\varepsilon(t)+f'(u)\cos\varepsilon(t)]E_{2}\nonumber\\
&+&\frac{1}{2}[f(u)-uf'(u)]E_{3}.
\end{eqnarray*}
Like we did before, from \eqref{3} and \eqref{dan} it is easy to conclude that
\begin{eqnarray}\label{eq2}
\varepsilon'(t)[1+f'(u)^{2}][u^{2}+f(u)^{2}]'=f''(u).
\end{eqnarray}

 Finally, we get the proof of Theorem \ref{t1} done.
\hfill $\Box$

 \vspace{.2in}

 \n \textbf{Proof of Theorem \ref{t2}:}
 
{\bf First Case:} We start this one by pointing out that, from \eqref{geohorizontal} we have that the horizontal straight lines are geodesics in $\mathcal{H}$. Now, from straightforward computation we can see that any base curve that is orthogonal to those horizontal straight lines (geodesics) has the following form
 \begin{eqnarray*}
 \alpha(u)=(-Bf(u),\,Cf(u),\,g(u)),
 \end{eqnarray*}
 where $f$ and $g$ are smooth functions.
 
 So, to begin with, consider a ruled surface in $\mathcal{H}$ given by 
 \begin{eqnarray}\label{regradaH}
 \Phi(u,\,v)=(-Bf(u),\,Cf(u),\,g(u))+v(C,\,B,\,0);\quad C^{2}+B^{2}=1.
 \end{eqnarray}
From \eqref{4}, considering the flow generated by $X_{2}$, we get
 \begin{eqnarray*}
 \Phi^{t}(u,\,v)=\left(\varepsilon(t)+Cv-Bf(u),\,Cf(u)+Bv,\,\frac{1}{2}\varepsilon(t)(Cf(u)+Bv)+g(u)\right).
 \end{eqnarray*}
 The first derivatives are 
 \begin{eqnarray}\label{d1}
  \Phi^{t}_{u}=-Bf'(u)E_{1}+Cf'(u)E_{2}
+\frac{1}{2}\left[2g'(u)-vf'(u)\right]E_{3}  
 \end{eqnarray}
 and
  \begin{eqnarray}\label{d2}
  \Phi^{t}_{v}= CE_{1}+BE_{2}+\frac{1}{2}f(u)E_{3}.
 \end{eqnarray}
We can conclude that
\begin{eqnarray*}
	E^{t}&=&\langle \Phi^{t}_{u},\,\Phi^{t}_{u}\rangle=f'(u)^{2}+\frac{1}{4}\left[2g'(u)-vf'(u)\right]^{2},\nonumber\\
	F^{t}&=&\langle \Phi^{t}_{u},\,\Phi^{t}_{v}\rangle=\frac{1}{4}f(u)[2g'(u)-vf'(u)],\nonumber\\
	G^{t}&=&\langle \Phi^{t}_{v},\,\Phi^{t}_{v}\rangle=1+\frac{1}{4}f(u)^{2}.
\end{eqnarray*}
Thus,  
\begin{eqnarray*}
E^{t}G^{t}-(F^{t})^{2}=E^{t}+(G^{t}_{u})^{2}.
\end{eqnarray*}
Let's compute the normal vector field $N^{t}$. From \eqref{basevetorial}, \eqref{d1} and \eqref{d2} we have
\begin{eqnarray}\label{ddd1}
N^{t}&=&\frac{1}{2}\|E^{t}G^{t}-(F^{t})^{2}\|^{-1/2}\Big\{\left[Cf(u)f'(u)-B(2g'(u)-vf'(u))\right]E_{1}\nonumber\\
&+&\left[Bf(u)f'(u)+C(2g'(u)-vf'(u))\right]E_{2}-2f'(u)E_{3}\Big\}.
\end{eqnarray}
 The coefficients of the second fundamental form are:
\begin{eqnarray*}
l^{t}&=&\langle \nabla_{\Phi^{t}_{u}}\Phi^{t}_{u},\,N^{t}\rangle=\frac{[f''(u)+\frac{1}{2}f(u)f'(u)^{2}][2g'(u)-vf'(u)]-f'(u)[2g''(u)-vf''(u)]}{2\|E^{t}G^{t}-(F^{t})^{2}\|^{1/2}}\nonumber\\
&=&\left\{F^{t}f'(u)^{2}+[g'(u)f''(u)-g''(u)f'(u)]\right\}\|E^{t}G^{t}-(F^{t})^{2}\|^{-1/2},\nonumber\\\\
m^{t}&=&\langle \nabla_{\Phi^{t}_{u}}\Phi^{t}_{v},\,N^{t}\rangle=\frac{[f(u)f'(u)]^{2}-[2g'(u)-vf'(u)]^{2}}{8\|E^{t}G^{t}-(F^{t})^{2}\|^{1/2}}\nonumber\\
&=&\frac{1}{2}(f'(u)^{2}G^{t}-E^{t})\|E^{t}G^{t}-(F^{t})^{2}\|^{-1/2},\nonumber\\\\
n^{t}&=&\langle \nabla_{\Phi^{t}_{v}}\Phi^{t}_{v},\,N^{t}\rangle=-\frac{f(u)[2g'(u)-vf'(u)]}{4\|E^{t}G^{t}-(F^{t})^{2}\|^{1/2}}
=-F^{t}\|E^{t}G^{t}-(F^{t})^{2}\|^{-1/2}.
\end{eqnarray*}
All the ingredients to provide the mean curvature \eqref{meancurva} were given, i.e.,
\begin{eqnarray}\label{m1}
H^{t}=\frac{\left[1+\frac{1}{4}f(u)^{2}\right][f''(u)g'(u)-f'(u)g''(u)]}{2\|E^{t}G^{t}-(F^{t})^{2}\|^{1/2}[E^{t}+(G^{t}_{u})^{2}]}.
\end{eqnarray}

On the other hand, 
\begin{eqnarray*}
 \partial_{t}\Phi^{t}=\varepsilon'(t)E_{1}+\varepsilon'(t)[Cf(u)+Bv]E_{3}.
\end{eqnarray*}
Therefore, 
\begin{eqnarray}\label{m2}
	\langle\partial_{t}\Phi^{t},\,N^{t}\rangle
	=\frac{-\varepsilon'(t)}{2\|E^{t}G^{t}-(F^{t})^{2}\|^{1/2}}[f'(u)(Af(u)+Bv)+2Bg'(u)].
\end{eqnarray}

Thus, combining \eqref{m1} and \eqref{m2} we obtain
\begin{eqnarray*}
-\varepsilon'(t)[f'(u)(Cf(u)+Bv)+2Bg'(u)]
=\frac{\left[4+f(u)^{2}\right][f''(u)g'(u)-f'(u)g''(u)]}{f'(u)^{2}[4+f(u)^{2}]+(2g'(u)-vf'(u))^{2}}.
\end{eqnarray*}
This proves the second item.

We pointed out that, the coefficients of the first and second fundamental forms for each item of this theorem are equal.

{\bf Second Case:} Therefore, to prove the ODE equation for $\Psi_{3}$, we just need to infer the left-hand side of \eqref{3}. So, from \eqref{4} and \eqref{regradaH} we get
 \begin{eqnarray*}
	\Phi^{t}(u,\,v)=(Cv-Bf(u),\,\varepsilon(t)+Bv+Cf(u),\,-\frac{1}{2}\varepsilon(t)[Cv-Bf(u)]+g(u)),
\end{eqnarray*}
and so
\begin{eqnarray*}
 \partial_{t}\Phi^{t}=\varepsilon'(t)\left[E_{2}+(Bf(u)-Cv)E_{3}\right].
\end{eqnarray*}
From the above identity and \eqref{ddd1} we have
\begin{eqnarray*}
	\langle\partial_{t}\Phi^{t},\,N^{t}\rangle=\frac{-\varepsilon'(t)}{2\|E^{t}G^{t}-(F^{t})^{2}\|^{1/2}}\left[f'(u)(Bf(u)-Cv)-2Cg'(u)\right].
\end{eqnarray*}

This implies that the same ODE gives the family of MCF for $X_{2}$ and $X_{3}$, provided a proper choice of $C$ and $B$.

{\bf Third Case:} Now, considering $\Psi_{4}$ in \eqref{4}, from \eqref{regradaH} we have
 \begin{eqnarray*}
\Phi^{t}(u,\,v)=(Cv-Bf(u),\,Bv+Cf(u),\,g(u)+\varepsilon(t)),
\end{eqnarray*}
and so
\begin{eqnarray*}
    \partial_{t}\Phi^{t}=\varepsilon'(t)E_{3}.
\end{eqnarray*}
Thus, 
\begin{eqnarray}\label{ddd}
    \langle \partial_{t}\Phi^{t},\,N^{t}\rangle=\frac{-\varepsilon'(t)f'(u)}{\|E^{t}G^{t}-(F^{t})^{2}\|^{1/2}}.
\end{eqnarray}
And thus, combining \eqref{m1} with \eqref{ddd} we get the fourth item.

{\bf Fourth Case:} To finish this paper, we prove the MCF for the first item of this theorem. It is easy to see that combining \eqref{4} and \eqref{regradaH} we obtain the parametrization of this soliton. So by Theorem \ref{t2}-1 we have
\begin{eqnarray*}
    \Phi^{t}_{u}&=&-f'(u)[C\sin\varepsilon(t)+B\cos\varepsilon(t)]E_{1}+f'(u)[C\cos\varepsilon(t)-B\sin\varepsilon(t)]E_{2}\nonumber\\   &+&\frac{1}{2}[2g'(u)-vf'(u)]E_{3}
\end{eqnarray*}
and
\begin{eqnarray*}
    \Phi^{t}_{v}&=&[C\cos\varepsilon(t)-B\sin\varepsilon(t)]E_{1}+[C\sin\varepsilon(t)+B\cos\varepsilon(t)]E_{2}+\frac{f(u)}{2}E_{3}.
\end{eqnarray*}

Thus, it is a straightforward computation to see that the coefficients of the first and second fundamental forms remain the same as in the previous cases.
However, we now can realize that the normal vector is different from the other cases, and is given by
\begin{eqnarray*}
    N^{t}&=&\frac{1}{2}\|E^{t}G^{t}-(F^{t})^{2}\|^{-1/2}\nonumber\\
    &\times&\Big\{[f(u)f'(u)(C\cos\varepsilon(t)-B\sin\varepsilon(t))-(2g'(u)-vf'(u))(C\sin\varepsilon(t)+B\cos\varepsilon(t))]E_{1}\nonumber\\
    &+&[f(u)f'(u)(C\sin\varepsilon(t)+B\cos\varepsilon(t))+(2g'(u)-vf'(u))(C\cos\varepsilon(t)-B\sin\varepsilon(t))]E_{2}\nonumber\\
    &-&2f'(u)E_{3}\Big\}.
\end{eqnarray*}

Moreover,
\begin{eqnarray*}
    \partial_{t}\Phi^{t}
    =\varepsilon'(t)\Big\{&-&[v(C\sin\varepsilon(t)+B\cos\varepsilon(t))+f(u)(C\cos\varepsilon(t)-B\sin\varepsilon(t))]E_{1}\nonumber\\
    &+&[v(C\cos\varepsilon(t)-B\sin\varepsilon(t))-f(u)(C\sin\varepsilon(t)+B\cos\varepsilon(t))]E_{2}\nonumber\\
    &-&\frac{1}{2}[v^{2}+f(u)^{2}]E_{3}\Big\}.
\end{eqnarray*}

%%%%%%%%%%%%%%%%%%%%%%%%%%%%%%%%%%%%%%%%%%%%%%
\iffalse
and
\begin{eqnarray*}
    \nabla_{\Phi_{v}^{t}}\Phi_{v}^{t}=\frac{f(u)}{2}\left[(A\sin\varepsilon(t)+B\cos\varepsilon(t))E_{1}+(B\sin\varepsilon(t)-A\cos\varepsilon(t))E_{2}\right].
\end{eqnarray*}
\fi
%%%%%%%%%%%%%%%%%%%%%%%%%%%%%%%%%%%%%%%%%%%%%%%

Then, from the last two idenetities we can calculate the left-hand side of  \eqref{3}, i.e.,
\begin{eqnarray*}
  \langle \partial_{t}\Phi^{t},\,N^{t}\rangle&=&  \frac{vg'(u)\varepsilon'(t)}{\|E^{t}G^{t}-(F^{t})^{2}\|^{1/2}}.
\end{eqnarray*}

Finally, combining \eqref{m1} with the above equation we have our result. 

 \hfill $\Box$
  \vspace{.2in}
\iffalse
\n \textbf{Acknowledgement:}
{ \it  The authors want to thanks the referee for his careful reading and helpful suggestions.}
\fi

\end{document}